\documentclass[11pt]{amsart}
\usepackage{amsthm,amsmath,amssymb,amscd,epsfig}

{\end{list}}
{
   \newtheorem{theorem}{Theorem}[section]
   \newtheorem{proposition}[theorem]{Proposition}
   \newtheorem{lemma}[theorem]{Lemma}

   \newtheorem{corollary}[theorem]{Corollary}

}
{\theoremstyle{definition}

}
{\theoremstyle{remark}

}

\newcommand{\RR}{{\mathbb{R}}}

\newcommand{\ZZ}{{\mathbb{Z}}}

\newcommand{\cA}{{\mathcal A}}

\newcommand{\cE}{{\mathcal E}}
\newcommand{\cF}{{\mathcal F}}

\newcommand{\cL}{{\mathcal L}}

\newcommand{\res}{\operatorname{res}}

\newcommand{\Star}{\operatorname{Star}}

\newcommand{\isom}{\simeq}

\begin{document}
\title{Ehrhart analogue of the $h$-vector.}

\author{Kalle Karu}
\thanks{The author was supported by an NSERC grant}
\address{Department of Mathematics\\ University of British Columbia \\
  1984 Mathematics Road\\
Vancouver, B.C. Canada V6T 1Z2}
\email{karu@math.ubc.ca}

\begin{abstract}
We consider a formula of Stanley that expresses the Ehrhart generating
polynomial of a polyhedral complex in terms of the $h$-polynomials of
toric varieties. We prove that the coefficients in this expression are
all non-negative and show that these coefficients can be found using the
decomposition theorem in intersection cohomology.
\end{abstract}

\maketitle

\section{Introduction}
\setcounter{equation}{0}

The title of this note comes from Example~7.13 in Stanley's paper
\cite{S}. In that example Stanley considers the Ehrhart problem of
counting lattice points in a polyhedral complex $C$. In analogy
with decomposing the $h$-vector of a subdivision of $C$ into
``local $h$-vectors'', he decomposes the Ehrhart generating polynomial
of $C$ into the same local $h$-vectors with coefficients
$c_{\sigma,i}$ for $\sigma\in C$, $i\geq 0$. We prove
Conjecture~7.14 in \cite{S} that these coefficients   $c_{\sigma,i}$
are non-negative by interpreting them in terms of orbifold cohomology.

Even though Stanley's definition of the numbers $c_{\sigma,i}$ is
combinatorial, it is more illuminating to define them using orbifold
cohomology. It is well-known that counting lattice points in a
polyhedral complex is equivalent to studying orbifold cohomology of
some toric orbifold, a relationship that is analogous to the equivalence
between counting faces of a (simplicial) polyhedral complex and
ordinary cohomology. In relating the combinatorics and cohomology we
follow the article of Mustata and Payne 
\cite{MP} where this equivalence is used to study a conjecture of Hibi.

Let us start by recalling the definition of orbifold cohomology
defined by Chen and Ruan \cite{CR}. Later we will specialize to the
case of toric varieties studied by Borisov, Chen and Smith \cite{BCS}.
 We are only interested in the dimensions of the orbifold
 cohomology spaces, not the ring structure. These dimensions were computed in
 a more general setting by Batyrev \cite{B} and
 Batyrev-Dais \cite{BD}.

Given a complete Gorenstein orbifold $X$, one decomposes its inertia
scheme:
\[ I(X) = \coprod_i X_i \] 
and defines the orbifold cohomology of $X$:
\[ H_{orb}^k(X;\RR) = \oplus_i H^{k-2s_i} (X_i;\RR),\]
where $s_i\geq 0$ is the ``age'' of the component $X_i$. If $Y\to X$
  is a crepant resolution of $X$, then the orbifold cohomology of $X$
  is isomorphic to the ordinary cohomology of $Y$. 

When $X_\Delta$ is a toric orbifold defined by a complete simplicial
Gorenstein fan $\Delta$, then $I(X_\Delta)$ is a union of orbit closures
$V_\sigma$ for 
$\sigma\in\Delta$ and we have 
\begin{equation} \label{def-orb}  H_{orb}^k(X_\Delta;\RR) =
  \bigoplus_{\sigma\in\Delta} \bigoplus_{i\geq 0} 
    [H^{k-2i} (V_\sigma;\RR)]^{c_{\sigma,i}}.
\end{equation}
The coefficients $c_{\sigma,i}$ have a combinatorial meaning: they
count lattice points in the interior 
of $Box(\sigma)$ (see Section~\ref{sec-fans} below for
the definition of $Box(\sigma)$).

Let now $X$ be an arbitrary complete Gorenstein variety that admits a
crepant resolution $Y\to X$, where $Y$ is an orbifold. Define the
{\em stringy cohomology} of $X$ by
\[ H_{str}^k(X;\RR) := H_{orb}^k(Y;\RR).\]
By a result of Yasuda \cite{Y} this does not depend on the choice of
$Y$. The dimensions of the stringy cohomology spaces 
$H_{str}^k(X;\RR)$ are given by Batyrev's stringy Betti numbers
\cite{B}. An arbitrary Gorenstein variety may not have any crepant orbifold
resolution, however, a Gorenstein toric variety $X_\Delta$ always has
one, hence the stringy cohomology $H_{str}^k(X_\Delta;\RR)$ is
well-defined. 

Note that orbifold cohomology was defined using the decomposition of
the inertia scheme. The behavior of the cohomology under crepant
morphisms was a theorem. For stringy cohomology we reverse the
situation: we define it using a crepant orbifold resolution. Its
decomposition then becomes a theorem.  

\begin{theorem}\label{thm-main} Let $X_\Delta$ be a complete Gorenstein
  toric variety. There exist non-negative integers $c_{\sigma,i}$
  such that 
\[  H_{str}^k(X_\Delta;\RR)  = \bigoplus_{\sigma\in\Delta} \bigoplus_{i\geq 0} 
    [IH^{k-2i} (V_\sigma;\RR)]^{c_{\sigma,i}},\qquad k\geq 0,\]
where $IH^* (V_\sigma)$ is the intersection cohomology of
$V_\sigma$. The numbers $c_{\sigma,i}$ depend on the cone $\sigma$
only, not on the fan $\Delta$.
\end{theorem}

The numbers  $c_{\sigma,i}$ are the ones calculated by Stanley in
\cite{S}. The theorem is proved by applying the decomposition
theorem in intersection cohomology to the crepant orbifold resolution
$Y\to X_\Delta$.

Let us now explain the relation with lattice point counting. If $P$ is
a $d$-dimensional lattice polytope in $\RR^n$, its Ehrhart generating
function is 
\[ F_P(t) = \sum_{j=0}^\infty \#(jP\cap \ZZ^n) t^j = \frac{\delta_0+\delta_1 t+
  \ldots + \delta_d t^d}{(1-t)^{d+1}}\]
for some non-negative integers $\delta_0,\ldots,\delta_d$. Similarly, when $C$
is a complex of lattice polytopes in $\RR^n$, such that $C$ is
topologically a $d$-sphere, then
\begin{equation}\label{eq-ehr}    F_C(t) = \sum_{j=0}^\infty \#(jC\cap
  \ZZ^n) t^j = \frac{\delta_0+\delta_1 t+ \ldots + \delta_{d+1}
  t^{d+1}}{(1-t)^{d+1}},
\end{equation}
with $\delta_0,\ldots,\delta_{d+1}$ non-negative integers. Denote
by $\delta_C(t)$ the polynomial in the numerator of the fraction in
(\ref{eq-ehr}).

We will consider the case where the complex $C$ is defined by a
complete Gorenstein fan $\Delta$ in $\RR^n$. Recall that $\Delta$ being
Gorenstein means that there exists a conewise linear integral function
$K_\Delta: \RR^n \to \RR$, such that $K_\Delta(v) = 1$ for all primitive
generators $v$ of the rays of $\Delta$. Let $C=K_\Delta^{-1}(1)$. Then
in the simplicial case the numbers $\delta_k$ in Equation~(\ref{eq-ehr})
are (see \cite{BCS, MP}) 
\[ \delta_k = \dim H_{orb}^{2k}(X_\Delta;\RR).\]
Since stringy cohomology was defined by a crepant resolution, for a
general $X_\Delta$ we have 
\[ \delta_k = \dim H_{str}^{2k}(X_\Delta;\RR).\]

For a toric variety $X$ write its $h$-polynomial:
\[ h_X(t) = \sum_{k\geq 0} \dim IH^{2k}(X;\RR) t^k.\]
Also let $c_\sigma(t)$ be the polynomial
\[  c_\sigma(t) = \sum_{j\geq 0} c_{\sigma,j} t^j,\]
with $c_{\sigma,j}$ defined in Theorem~\ref{thm-main}. Then
Theorem~\ref{thm-main} implies by computing dimensions on both sides:

\begin{corollary}\label{cor-main} \[ \delta_C(t) = \sum_{\sigma\in\Delta} c_\sigma(t)  h_{V_\sigma}(t).\] 
\end{corollary}

Note that in the corollary all polynomials have non-negative
coefficients. The polynomial $c_\sigma(t)$ depends on the cone
$\sigma$ only, while the polynomial $h_{V_\sigma}(t)$ depends on the
poset of $\Star_\Delta(\sigma)$ only.

To decompose the stringy cohomology as a direct sum of intersection
cohomologies, we use the combinatorial theory of locally free and
flabby sheaves on the fan $\Delta$ (\cite{BBFK, BL}). We construct a
locally free and flabby sheaf $\cE_\Delta$, called the Ehrhart sheaf, such
that the global sections of $\cE_\Delta$ give the equivariant stringy
cohomology of $X_\Delta$. Applying the combinatorial decomposition
theorem to $\cE_\Delta$ then proves Theorem~\ref{thm-main}.

Let us mention a few generalizations of Theorem~\ref{thm-main}. The
construction of the Ehrhart sheaf $\cE_\Delta$ and the decomposition
of this sheaf makes sense for any (rational) fan.
Theorem~\ref{thm-main} and Corollary~\ref{cor-main} hold, for example,
if the fan $\Delta$ is quasi-convex \cite{BBFK}. With the same tools
one can even treat the case of abstract fans corresponding to abstract
polyhedral complexes, but we will not pursue these generalizations here.
Let us only bring the analogue of Corollary~\ref{cor-main} in the
following quasi-convex case.

Suppose the complex $C$ consists of a single lattice polytope
$P$. Then the fan $\Delta$ consists of a single cone over the polytope
$P$ and all its faces. The $h$-polynomial of the affine toric variety
$X_\Delta$ is usually  called the $g$-polynomial of the polar dual polytope
$P^\circ$. Similarly, the $h$-polynomial of $V_\sigma$ for a cone $\sigma$
corresponding to a face $F\leq P$ is the $g$-polynomial of the dual face
$F^* \leq P^\circ$. The formula in  Corollary~\ref{cor-main} can now be written
as:
\[ \delta_P(t) = \sum_{F \leq P} c_F(t) g_{F^*}(t),\]
where we replaced  $c_\sigma(t)$ by the corresponding $c_F(t)$. Again,
all polynomials in this formula have non-negative integer
coefficients. The polynomial $c_F(t)$ depends on the face $F$
only, the polynomial $g_{F^*}(t)$ depends on the poset of the dual
face $F^*$ only.

\section{Combinatorial decomposition theorem}

We recall briefly the decomposition theorem for locally free and
flabby sheaves on a fan $\Delta$  (\cite{BBFK, BL}). All vector spaces
are over the field $\RR$.

Let $\Delta$ be a polyhedral fan in $\RR^n$. The fan (as a finite set
of cones) is given the topology in which open sets are subfans of
$\Delta$. A sheaf of vector spaces $F$ in this topology consists of 
the data:
\begin{itemize}
\item A vector space $F_\sigma$ for each $\sigma\in\Delta$, the
  stalk of $F$ at $\sigma$.
\item A linear map $res^\sigma_\tau: F_\sigma\to F_\tau$ for $\tau$ a
  face of $\sigma$, such that $res^\sigma_\sigma = Id$ and
  $\res^\tau_\rho \circ res^\sigma_\tau = res^\sigma_\rho$ whenever
  $\sigma>\tau>\rho$. These maps are called restriction maps.
\end{itemize}

To a sheaf $F$ on $\Delta$ we can apply the usual
constructions, such as taking global sections or computing the sheaf
cohomology. Given a subdivision of fans $\phi: \hat{\Delta}\to \Delta$ and a 
sheaf $F$ on $\hat{\Delta}$, we let $\phi_*(F)$ be the push-forward of
$F$. 

Recall that a sheaf $F$ is flabby if for any open sets $U\subset V$,
sections on $U$ can be lifted to sections on $V$. In the case of fans
this amounts to the surjectivity of the maps
\[ F_\sigma \to \Gamma(F, \partial\sigma) \] 
for all cones $\sigma\in\Delta$.

Let $\cA = \cA_\Delta$ be the sheaf defined by:
\begin{itemize}
\item $\cA_\sigma$  is the space of polynomial functions on $\sigma$.
\item $res^\sigma_\tau$ is the restriction of functions.
\end{itemize}
The sheaf $\cA$ is a graded sheaf of rings. It is flabby if and only
if the fan $\Delta$ is simplicial. If the fan is complete and
simplicial, then global sections of $\cA$ form the $T$-equivariant
cohomology ring of $X_\Delta$, where $T$ is the torus.

Let $\cF$ be a sheaf of $\cA$-modules. $\cF$ is called locally free if
$\cF_\sigma$ is a finitely generated free $\cA_\sigma$-module for all
$\sigma\in\Delta$. The decomposition theorem states that a locally
free flabby sheaf $\cF$ can be decomposed as a direct sum of
elementary sheaves:
\[ \cF = \bigoplus_{\sigma\in\Delta} \bigoplus_i
\cL_\sigma^{c_{\sigma,i}}[i].\]
The sheaves $\cL_\sigma$ are indecomposable locally free flabby sheaves
supported on $\Star \sigma = \{\pi\geq\sigma\}$. In the simplicial
case $\cL_\sigma$ is simply the sheaf $\cA$ restricted to  $\Star
\sigma$. The non-negative integers $c_{\sigma,i}$ depend on the
restriction of $\cF$ to the cone $\sigma$ and all its faces only.

Assume now that $\Delta$ is complete. Then the global sections of the
sheaf $\cL_\sigma$ form the $T$-equivariant intersection cohomology of the
orbit closure $V_\sigma$. The equivariant cohomology is related to the
non-equivariant cohomology as follows. Let $A$ be the ring of global
polynomial functions on $\RR^n$. Then $\cL_\sigma$ is a sheaf of
$A$-modules. The space of global sections of $\cL_\sigma$ forms a free
$A$-module and we have
\[ \Gamma(\cL_\sigma,\Delta) = IH_T^*(V_\sigma;\RR) \isom
IH^*(V_\sigma;\RR)\otimes_\RR A.\]
The degree convention is that a section of degree $k$ gives a
cohomology class of degree $2k$. In particular, the Hilbert function
of the space of global sections of $L_\sigma$ is
\[ \sum_{j=0}^\infty \dim \Gamma(\cL_\sigma,\Delta)_j t^j =
\frac{h_0+h_1 t+ \ldots + h_{n} t^{n}}{(1-t)^{n}},\]
where $h_k = \dim IH^{2k}(V_\sigma;\RR)$.

\section{Gorenstein fans} \label{sec-fans}

Let $\sigma$ be a rational polyhedral pointed cone in $\RR^n$. Let
$v_1,\ldots,v_m$ be the primitive lattice points on the edges of
$\sigma$. Recall that $\sigma$ is called Gorenstein if there exists a
linear function $K_\sigma: \sigma\to \RR$ taking integral values on
$\ZZ^n\cap \sigma$ and such that $K_\sigma(v_i) = 1$ for
$i=1,\ldots,m$. A fan $\Delta$ is Gorenstein if all cones
$\sigma\in\Delta$ are Gorenstein. On a Gorenstein fan the functions
$K_\sigma$ glue to a continuous conewise linear function $K_\Delta$.

Consider a Gorenstein $d$-dimensional cone $\sigma$ and the power series
\[ F_\sigma(t) = \sum_{j=0}^\infty \#(K_\sigma^{-1}(j)\cap \ZZ^n)
t^j.\]
Since this is the Ehrhart generating function for the polytope
$P=K_\sigma^{-1}(1)$, we can write it as a rational function
\[ F_\sigma(t) = \frac{\delta_0+\delta_1 t+ \ldots + \delta_{d-1} t^{d-1}}{(1-t)^{d}}.\]
Note that if $\cE_\sigma$ is the free graded $\cA_\sigma$-module 
\[ \cE_\sigma = \cA_\sigma^{\delta_0} \oplus \cA_\sigma^{\delta_1}[-1] \oplus
\cdots \cA_\sigma^{\delta_{d-1}}[-(d-1)],\]
then its Hilbert series is precisely $F_\sigma(t)$. We construct a
locally free and flabby sheaf $\cE$ on $\Delta$ with stalks
$\cE_\sigma$ as above and call it the Ehrhart sheaf. By the
decomposition theorem there is up to an isomorphism a unique locally
free flabby sheaf on $\Delta$ with the given stalks, hence the Ehrhart
sheaf is unique.

Let us start with the case where $\Delta$ is
simplicial. Let $\sigma$ be a simplicial cone with primitive
generators $v_1,\ldots,v_d$, and let $Box(\sigma)$ be
\[ Box(\sigma) = \{v\in\ZZ^n | v = \alpha_1 v_1+ \cdots +\alpha_d v_d
\text{ for some } 0\leq \alpha < 1\}.\]
Define $\cE_\sigma$ to be the free graded $\cA_\sigma$-module with basis
$Box(\sigma)$, where a  basis element $v\in Box(\sigma)$ has degree
$K_\sigma(v)$. (This definition of  $\cE_\sigma$ agrees with the one
given above.) If $\tau$ is a face of $\sigma$ then $Box(\tau)\subset
Box(\sigma)$, and we get the restriction map $\cE_\sigma\to \cE_\tau$
by sending  a basis element $v\in Box(\sigma)$ to $v$ if  $v\in
Box(\tau)$ and to zero otherwise. This defines a locally free sheaf
$\cE = \cE_\Delta$ on $\Delta$.

\begin{lemma} The sheaf $\cE_\Delta$ is flabby on the simplicial fan
  $\Delta$.
\end{lemma}

\begin{proof}
Note that $\cE_\Delta$ decomposes into a finite direct sum
\[ \cE_\Delta = \bigoplus_{v\in\ZZ^n} \cE_v,\]
where 
\[ \cE_{v,\sigma} = \begin{cases} \cA_\sigma & \text{ if $v\in
    Box(\sigma)$,} \\
0 & \text{ otherwise.}
\end{cases} \]
If $\sigma$ is the smallest cone such that $v\in Box(\sigma)$, then
\[ \cE_v = \cA|_{\Star \sigma} = \cL_\sigma,\]
which is flabby.
\end{proof}

We note that the decomposition $\cE = \oplus_{v} \cE_v$ in the proof
above corresponds to the decomposition (\ref{def-orb}) of the orbifold
cohomology of $X_\Delta$. 

Let $\Gamma(\cE_\Delta,\Delta)_k$ be the degree $k$ component of the
graded vector space of global sections.

\begin{lemma} We have 
\[ \dim   \Gamma(\cE_\Delta,\Delta)_k = \#(K_\Delta^{-1}(k)\cap
\ZZ^n).\]
\end{lemma}

\begin{proof} Define a sheaf $F$ of graded vector spaces on
  $\Delta$ as follows. The stalk $F_\sigma$ has basis $\sigma\cap\ZZ^n$ with
  $v\in \sigma\cap\ZZ^n$ having degree $K_\sigma(v)$. The restriction
  map $res^\sigma_\tau$ sends a basis element $v$ to the same basis
  element or to zero as appropriate. Clearly $F$ is a flabby sheaf of
  vector spaces on $\Delta$ satisfying the equality stated for
  $\cE_\Delta$. Similarly to the case of locally free flabby $\cA$-modules,
  there is up to an isomorphism a unique graded flabby sheaf of vector
  spaces on $\Delta$ with the given stalks. Hence it suffices to prove
  that $\cE_\Delta$ and $F$ have the same stalks. This is clear. 
\end{proof}

Let now $\Delta$ be an arbitrary rational polyhedral complete
Gorenstein fan in $\RR^n$. Then $\Delta$ has a simplicial subdivision
\[ \phi: \hat\Delta\to \Delta,\]
such that $K_{\hat\Delta} = K_\Delta \circ\phi$. Such a subdivision
$\phi$ is called crepant. Define 
\[ \cE_\Delta = \phi_* \cE_{\hat\Delta}.\]

\begin{proposition} $\cE_\Delta$ is a locally free flabby sheaf on
  $\Delta$, such that 
\[ \dim\Gamma(\cE_\Delta,\Delta)_k = \#(K_\Delta^{-1}(k)\cap \ZZ^n).\]
\end{proposition}

\begin{proof} $\cE_\Delta$ is a locally free and flabby by \cite{BBFK,
    BL}. The second statement follows from the equality 
\[ \Gamma(\cE_\Delta,\Delta)_k =
\Gamma(\cE_{\hat\Delta},\hat\Delta)_k\] 
and the previous lemma.
\end{proof}

As explained in the introduction, the proposition states that
$\Gamma(\cE_\Delta,\Delta)$ gives the equivariant stringy cohomology
of $X_\Delta$. Since $\cE_\Delta$ is locally free and flabby, we
decompose it as
\[ \cE_\Delta = \bigoplus_{\sigma\in\Delta} \bigoplus_i
\cL_\sigma^{c_{\sigma,i}}[i], \]
for some non-negative integers $c_{\sigma,i}$ that depend on the cone
$\sigma$ only. Taking global sections of both sides proves
Theorem~\ref{thm-main}.


\begin{thebibliography}{CC}

\bibitem{BBFK} G.\ Barthel, J.-P.\ Brasselet, K.-H.\ Fieseler and L.\ Kaup,
{\em Combinatorial Intersection cohomology for Fans}, T\^ohoku
Math.\ J.\ {\bf 54} (2002) 1--41.

\bibitem{B} V. Batyrev, {\em Stringy Hodge numbers of varieties with
  Gorenstein canonical singularities}, in Integrable systems and
  algebraic geometry (Kobe/Kyoto, 1997), 1-32, World Sci. Publishing,
  River Edge, NJ, 1998.

\bibitem{BD} V. Batyrev and D. Dais, {\em Strong McKay correspondence,
  string-theoretic Hodge numbers and mirror symmetry}, Topology 35
  (1996), 901-929. 

\bibitem{BCS} L. Borisov, L. Chen and G. Smith, {\em The orbifold
  Chow ring of toric Deligne-Mumford stacks},
J. Amer. Math. Soc. 18 (2005), 193-215.

\bibitem{BL} P.\ Bressler and V.\ Lunts, {\em Intersection cohomology
  on nonrational polytopes},  Compositio Math.\  {\bf 135}  (2003),  no.\ 3,
  245-278.

\bibitem{CR} W. Chen and Y. Ruan, {\em A new cohomology theory of
  orbifold}, Comm. Math. Phys. 248 (2004), 1-31.

\bibitem{MP} M. Mustata and S. Payne, {\em Ehrhart polynomials and stringy Betti numbers}, math.AG/0504486. Math. Ann., to appear.

\bibitem{S} R. Stanley, {\em Subdivisions and local h-vectors},
  J. Amer. Math. Soc. 5 (1992), 805-851.

\bibitem{Y} T. Yasuda, {\em Twisted jets, motivic measures and
  orbifold cohomology}, Compos. Math. 140 (2004), 396-422.



\end{thebibliography}
\end{document}